\documentstyle[amsfonts,11pt]{article}

\newtheorem{theo}{Theorem}[section]

\newtheorem{remar}[theo]{Remark}
\newtheorem{prop}[theo]{Proposition}

\newtheorem{lemma}[theo]{Lemma}

\def\Ric{\mathop{\hbox{Ric}}}
\def\scal{\mathop{\hbox{scal}}}

\newcommand{\fdim}{\hspace*{\fill}$\Box$}
\newcommand{\dimostr}{{\bf Proof: }}

\newcommand{\R}{\Bbb{R}}

\newcommand{\complex}{\Bbb{C}}
\newcommand{\C}{\Bbb{C}}

\newcommand{\K }{K\"{a}hler}

\begin{document}

\centerline {\LARGE\bf Canonical metrics on Hartogs domains
\footnote{The authors were supported
by the M.I.U.R. Project \lq\lq Geometric Properties of Real and
Complex Manifolds''.}}

\vspace{0.5cm}

\centerline{\small Andrea Loi and Fabio Zuddas } \centerline{\small Dipartimento di
Matematica e Informatica -- Universit\`{a} di Cagliari -- Italy}
\centerline{\small e-mail address: loi@unica.it, fzuddas@unica.it}

\vspace{0.3cm}

\begin{abstract}
An $n$-dimensional Hartogs domain $D_F$  with strongly
pseudoconvex boundary can be equipped with a natural \K\ metric
$g_F$. This paper contains two results. In the first one we prove
that if  $g_F$ is  an extremal \K\ metric then $(D_F, g_F)$ is
holomorphically isometric  to an open subset of the $n$-dimensional
complex hyperbolic space. In the second one we prove the same
assertion under the assumption that there exists a real
holomorphic vector field $X$ on $D_F$ such that $(g_F, X)$ is a
\K--Ricci soliton.

\K--Ricci solitons; generalized scalar curvatures; canonical metrics.

{\it{Subj.Class}}: 53C55, 32Q15, 32T15.

\end{abstract}

\section{Introduction and  statements of the main results}
The study of the existence and uniqueness of a  preferred \K\
metric on a given  complex manifold $M$ is a very interesting and
important  area of research, both from the mathematical and  from
the physical point of view. Many definitions  of    canonical
metrics  (Einstein, constant scalar curvature, extremal, \K--Ricci
solitons and so on) have been given both in the compact and in the
noncompact case (see e.g. \cite{besse},  \cite{gau} and
\cite{tian}). In the noncompact case many important questions are
still open. For example Yau raised the question on the
classification of Bergman Einstein metrics on strongly
pseudoconvex domains and S.--Y.  Cheng conjectured that if the
Bergman metric on a strongly pseudoconvex domain is Einstein, then
the domain is biholomorphic to the ball (see \cite{si}).

\vskip 0.3cm

In this paper we are interested in extremal \K\
metrics and \K--Ricci solitons  on  a  particular class  of strongly
pseudoconvex domains, the so called  {\em Hartogs domains} (see the next
section   for their definition and main properties). 

Our main results  are the following theorems.

\begin{theo}\label{mainteor}
Let $(D_F, g_F)$ be an $n$-dimensional strongly pseudoconvex Hartogs domain. Assume that
$g_F$ is an extremal \K\ metric.
Then $(D_F, g_F)$ is  holomorphically isometric to an open subset of the  $n$-dimensional complex hyperbolic space.
\end{theo}

\begin{theo}\label{mainteor2}
Let $(D_F, g_F)$ be an $n$-dimensional strongly pseudoconvex
Hartogs domain and let $X$ be a real holomorphic vector field on
$D_F$ such that $(g_F, X)$ is a \K--Ricci soliton. Then $g_F$ is
\K--Einstein. Consequently, $(D_F, g_F)$ is  holomorphically
isometric to an open subset of the $n$-dimensional complex
hyperbolic space.
\end{theo}

Notice that  (compare with  Cheng's conjecture above) the
assumptions on the metric $g_F$ in Theorem \ref{mainteor} and
Theorem \ref{mainteor2} are  weaker than Einstein's condition. To
this regard it is worth pointing out that when $g_F$ equals the
Bergman metric on $D_F$, then $(D_F, g_F)$ is holomorphically
isometric to an open subset of the complex hyperbolic space (see
Theorem 1.3 in \cite{riemhartogs} for a proof).

\vskip 0.3cm

The paper is organized as follows. In the next section, after
recalling the definition of Hartogs domains, we analyze their
pseudoconvexity,   and we  prove  a
lemma regarding their generalized scalar curvatures.
Sections \ref{final} and \ref{final2} are dedicated to the proofs
of Theorem \ref{mainteor} and Theorem \ref{mainteor2}
respectively.

\section{Strongly pseudoconvex Hartogs domains}\label{pseudo}
Let $x_0 \in \R^+ \cup \{ + \infty \}$ and let $F: [0, x_0)
\rightarrow (0, + \infty)$ be a decreasing continuous function,
smooth on $(0, x_0)$. The {\em Hartogs domain} $D_F\subset
{\C}^{n}$ associated to the function $F$ is defined by
$$D_F = \{ (z_0, z_1,\dots ,z_{n-1}) \in {\C}^{n} \; | \; |z_0|^2 < x_0, \ |z_1|^2 +\cdots+ |z_{n-1}|^2 < F(|z_0|^2)
\}.$$

One can prove that  the assumption of strongly pseudoconvexity of
$D_F$ is equivalent (see Proposition \ref{Kmet} below) to the fact
that the natural $(1, 1)$-form on $D_F$  given by
\begin{equation}\label{omegaf}
\omega_F = \frac{i}{2} \partial \overline{\partial}
\log \frac{1}{F(|z_0|^2) - |z_1|^2 - \cdots-|z_{n-1}|^2}
\end{equation}
is  a \K\ form on $D_F$.
The \K\  metric  $g_F$ associated to the \K\ form $\omega_F$
is the metric we will be  dealing with in the present paper.
Observe that for $F(x)=1-x, 0\leq x< 1$, $D_F$ equals the $n$-dimensional complex hyperbolic space
${\complex} H^n$ and $g_F$ is the hyperbolic metric, i.e. $g_F=g_{hyp}$.
In the $2$-dimensional case this metric has been considered in \cite{eng}
and \cite{constscal} in the framework of quantization of \K\ manifolds.
In  \cite{hartogs},  the first author studied the \K\ immersions of
$(D_F, g_F)$
into finite or infinite dimensional complex space forms, \cite{alcu}  is concerned with
the existence of global  symplectic coordinates on $(D_F, \omega_F)$ and
\cite{riemhartogs} deals with the Riemannian geometry of $(D_F, g_F)$
(in particular in this paper  one can find necessary and sufficient conditions in terms of $F$ for the completeness of
the metric $g_F$).
\begin{prop}\label{Kmet}
Let $D_F$ be a Hartogs domain in ${\complex}^{n}$. Then the
following conditions are equivalent:
\begin{itemize}
\item [(i)] the $(1, 1)$-form $\omega_F$  given by (\ref{omegaf})
is a \K\ form; \item [(ii)] the function $- \frac{x F'(x)}{F(x)}$
is strictly increasing, namely $-( \frac{x F'(x)}{F(x)})^{'} >0$
for every $x \in [0, x_0)$; \item [(iii)] the boundary of $D_F$ is
strongly pseudoconvex at all $z = (z_0, z_1,\dots,z_{n-1})$ with
$|z_0|^2 < x_0$ ; \item [(iv)] $D_F$ is strongly pseudoconvex.
\end{itemize}
\end{prop}
 \dimostr
 $(i)\Leftrightarrow (ii)$
 Set
 \begin{equation}\label{eqA}
 A = F(|z_0|^2)
- |z_1|^2 - \cdots - |z_{n-1}|^2.
\end{equation}
 Then $\omega_F$ is a \K\ form if and only if
the real-valued function $\Phi = -\log A$ is strictly
plurisubharmonic, i.e. the matrix $g_{\alpha\bar
\beta}=(\frac{\partial^2 \Phi}{\partial z_{\alpha}\partial \bar
z_{\beta}})$, $\alpha,\beta=0,\dots,n-~1$ is positive definite,
where
\begin{equation}\label{kform}
\omega_F=\frac{i}{2}\sum_{\alpha, \beta=0}^{n-1}g_{\alpha\bar \beta}dz_{\alpha}\wedge d\bar z_{\beta}.
\end{equation}
A straightforward computation gives
$$\frac{\partial^2 \Phi}{\partial z_0 \partial \bar z_0} =
\frac{F'^2(|z_0|^2) |z_0|^2 - ( F''(|z_0|^2) |z_0|^2 +
F'(|z_0|^2))A}{A^2} , $$
$$\frac{\partial^2 \Phi}{\partial z_0 \partial \bar z_\beta} =
- \frac{F'(|z_0|^2) \bar{z_0} z_{\beta}}{A^2}, \ \ \beta=1,\dots,n-1$$
and
$$\frac{\partial^2 \Phi}{\partial z_{\alpha} \partial \bar z_{\beta}} =
\frac{\delta_{\alpha \beta} A + \bar{z_{\alpha}} z_{\beta}}{A^2},
\ \ \alpha,\beta=1,\dots,n-1 .$$

Then, by setting
\begin{equation}\label{eqc}
C = F'^2(|z_0|^2) |z_0|^2 - ( F''(|z_0|^2) |z_0|^2 +
F'(|z_0|^2))A,
\end{equation}
one sees that the matrix $h=(g_{\alpha\bar
\beta})=(\frac{\partial^2 \Phi}{\partial z_{\alpha} \partial \bar
z_{\beta}})_{\alpha, \beta=0,\dots ,n-1}$  is given by:

\begin{small}
\begin{equation}\label{matrixg}
h= \frac{1}{A^2} \left( \begin{array}{cccccc}
C & - F' \bar z_0 z_1 & \dots  & - F' \bar z_0 z_{\alpha}&\dots & - F' \bar z_0 z_{n-1}\\
- F' z_0 \bar z_1  & A + |z_1|^2 &\dots  & \bar z_1 z_{\alpha} &\dots& \bar z_1 z_{n-1}\\
\vdots&\vdots &&\vdots&&\vdots \\
-F' z_0 \bar z_{\alpha} & z_1 \bar z_{\alpha} &\dots & A +
|z_{\alpha}|^2 &\dots&
 \bar z_{\alpha}z_{n-1}\\
\vdots&\vdots &&\vdots&&\vdots \\
-F' z_0 \bar z_{n-1} & z_1 \bar z_{n-1} &\dots & z_{\alpha} \bar
z_{n-1} &\dots & A + |z_{n-1}|^2
\end{array} \right) .
\end{equation}
\end{small}

First notice that the $(n-1)\times (n-1)$  matrix obtained by
deleting the first row and the first  column of $h$ is positive
definite. Indeed it is not hard to see that, for all $1\leq
\alpha\leq n-1$,
$$\det \left( \begin{array}{cccc} A + |z_{\alpha}|^2 & \bar
z_{\alpha} z_{\alpha+1} &\dots &
\bar z_{\alpha} z_{n-1}\\
 \vdots&\vdots &&\vdots \\
\bar z_{n-1} z_{\alpha} & \bar z_{n-1} z_{\alpha+1} &\dots& A +
|z_{n-1}|^2
\end{array} \right)=$$
\begin{equation}\label{minori principali}
 =A^{n- \alpha} + A^{n- \alpha
-1}(|z_{\alpha}|^2+\cdots+|z_{n-1}|^2) >0 .
\end{equation}

On the other hand, by the  Laplace expansion along the first row,
we get

$$\det(h) = \frac{C}{A^{2n}} [A^{n-1} + A^{n-2}(|z_1|^2+\cdots+|z_{n-1}|^2)] + $$

$$ + \frac{F' \bar z_0 z_1}{A^{2n}} \det \left( \begin{array}{cccc} -F' z_0 \bar z_1 & z_2 \bar z_1 &\dots&
z_{n-1} \bar z_{1}\\ -F' z_0 \bar z_2 & A + |z_2|^2 &\dots & z_{n-1} \bar
z_2\\
\vdots&\vdots &&\vdots \\
 -F' z_0 \bar z_{n-1} & z_2 \bar z_{n-1} &\dots & A + |z_{n-1}|^2
\end{array} \right) +\cdots + $$

$$+ (-1)^{n} \frac{F' \bar z_0 z_{n-1}}{A^{2n}} \det \left( \begin{array}{cccc} -F' z_0 \bar z_1 & A + |z_1|^2 &\dots&
z_{n-2} \bar z_{1}\\ -F' z_0 \bar z_2 & z_1 \bar z_2 &\dots& z_{n-2} \bar z_2\\
\vdots&\vdots &&\vdots \\
-F' z_0 \bar z_{n-1} & z_1 \bar z_{n-1} &\dots & z_{n-2} \bar
z_{n-1}
\end{array} \right) = $$

$$ = \frac{C}{A^{2n}} [A^{n-1} + A^{n-2}(|z_1|^2+\cdots +|z_{n-1}|^2)] + $$

$$ + \frac{F'^2 |z_0|^2 |z_1|^2}{A^{2n}} \det \left( \begin{array}{cccc} - 1 & z_2 &\dots&
z_{n-1}\\ - \bar z_2 & A + |z_2|^2 &\dots& z_{n-1} \bar z_2\\
\vdots&\vdots &&\vdots \\
- \bar z_{n-1} & z_2 \bar z_{n-1} &\dots& A + |z_{n-1}|^2
\end{array} \right) +\cdots+ $$

$$+ (-1)^{n} \frac{F'^2 |z_0|^2 |z_{n-1}|^2}{A^{2n}}
\det \left( \begin{array}{cccc} - \bar z_1 & A + |z_1|^2 &\dots&
z_{n-2} \bar z_{1}\\ - \bar z_2 & z_1 \bar z_2 &\dots& z_{n-2} \bar z_2\\
\vdots&\vdots &&\vdots \\
- 1 & z_1 &\dots& z_{n-2}
\end{array} \right) = $$

$$\frac{1}{A^{n+2}}[ C A + (C - F'^2 |z_0|^2 )(
|z_1|^2+\cdots +|z_{n-1}|^2)].$$

By substituting  (\ref{eqA}) and (\ref{eqc}) into this last
equality one gets
\begin{equation}\label{determinante metrica}
\det(h) = - \frac{F^2}{A^{n+1}} \left( \frac{x F'}{F}
\right)'|_{x=|z_0|^2} .
\end{equation}

Hence, by (\ref{minori principali}) and (\ref{determinante
metrica}),  the matrix $(\frac{\partial^2 \Phi}{\partial z_{\alpha}
\partial \bar z_{\beta}})$ is positive definite if and only if
$\left( \frac{x F'}{F} \right)'<0$.
\vskip 0.3cm

Before proving equivalence
$(ii)\Leftrightarrow (iii)$ we briefly recall some facts on complex domains
(see e.g. \cite{gurossi}).
Let  $\Omega \subseteq {\C}^{n}$
 be any complex domain of $ {\C}^{n}$   with smooth boundary $\partial\Omega$, and let $z \in \partial \Omega$.
 Assume that  there exists a smooth function
 $\rho : {\C}^n \rightarrow
{\R}$  (called  {\em defining function for $\Omega$ at $z$})
satisfying the following: for some neighbourhood $U$ of $z$, $\rho
< 0$ on $U \cap \Omega$, $\rho
>0$ on $U \setminus \overline \Omega$ and $\rho = 0$ on
$U \cap \partial \Omega$; $grad \ \rho \neq 0 $ on $\partial
\Omega$. In this case $\partial \Omega$ is said to be {\em
strongly pseudoconvex at $z$} if the {\it Levi form}

$$L(\rho, z)(X) = \sum_{\alpha,\beta=0}^{n-1} \frac{\partial^2 \rho}{\partial z_{\alpha} \partial \bar z_{\beta}}(z) X_{\alpha}\bar X_{\beta} $$
is positive definite on
$$S_{\rho}=\{ (X_0,\dots,X_{n-1}) \in {\C}^{n} \ |
\ \sum_{\alpha=0}^{n-1} \frac{\partial \rho}{\partial z_{\alpha}}(z) X_{\alpha} = 0 \}$$
(it is easily seen that this definition does not depend on the particular  defining function $\rho$).

\vskip 0.1cm

\noindent
$(ii)\Leftrightarrow (iii)$ Let now  $\Omega = D_F$ and
let us  fix   $z = (z_0, z_1,\dots,z_{n-1}) \in \partial D_F$ with  $|z_0|^2 < x_0$.
Then,   $|z_1|^2
+\cdots +|z_{n-1}|^2 = F(|z_0|^2)$.
In this case
$$\rho(z_0,z_1,\dots ,z_{n-1}) = |z_1|^2 +\cdots +|z_{n-1}|^2 - F(|z_0|^2)$$
 is a (globally)
defining function for $D_F$ at $z$, the Levi form for $D_F$ reads
as
\begin{equation}\label{levi1}
L(\rho, z)(X) = |X_1|^2 + \cdots + |X_{n-1}|^2 - (F' + F'' |z_0|^2)
|X_0|^2
\end{equation}
and
\begin{equation}\label{levi2}
S_{\rho}=\{ (X_0, X_1,\dots ,X_{n-1}) \in {\C}^{n} \ | -F' \bar z_0 X_0 +
\bar z_1 X_1 +\cdots + \bar z_{n-1} X_{n-1} = 0 \}.
\end{equation}

We distinguish two cases: $z_0=0$ and $z_0\neq 0$.
At  $z_0 = 0$ the Levi form reads as
$$L(\rho, z)(X) = |X_1|^2 +\cdots + |X_{n-1}|^2 - F'(0) |X_0|^2 $$
which is strictly positive for {\em any}  non-zero vector $(X_0,
X_1,\dots , X_{n-1})$ (not necessarily in $S_{\rho}$)
 because $F$ is assumed to be decreasing.

 If $z_0\neq 0$  by (\ref{levi2}) we obtain
 $X_0=  \frac{\bar z_1 X_1 +\cdots + \bar z_{n-1} X_{n-1}}
 {F' \bar z_0}$
 which,  substituted in (\ref{levi1}),  gives:
\begin{equation}\label{Levi ristretto}
L(X, z) = |X_1|^2 +\cdots + |X_{n-1}|^2 - \frac{F' + F'' |z_0|^2}{F'^2
|z_0|^2} |\bar z_1 X_1 +\cdots + \bar z_{n-1} X_{n-1} |^2.
\end{equation}

Therefore we are reduced  to show that:

\vskip 0.1cm

{\em
$(xF' / F)' <0$ for $x \in (0, x_0)$
if and only if
$L(X, z)$ is strictly positive
for every $(X_1,\dots ,X_{n-1}) \neq (0,\dots
, 0)$ and every $(z_0,z_1,\dots,z_{n-1}) \in \partial D_F$, $0 <
|z_0|^2 < x_0$.}

\vskip 0.3cm

If $(x F' / F)' <0$ then $(F' + x F'')F < x F'^2$ and,  since
$F(|z_0|^2) = |z_1|^2 +\cdots + |z_{n-1}|^2$, we get:
$$L(X, z) > |X_1|^2 +\cdots +
|X_{n-1}|^2 - \frac{1}{F(|z_0|^2)} |\bar z_1 X_1 +\cdots + \bar z_{n-1} X_{n-1}
|^2=$$
$$= \frac{(|X_1|^2 + \cdots + |X_{n-1}|^2)(|z_1|^2 + \cdots  + |z_{n-1}|^2) - |\bar z_1 X_1
+\cdots + \bar z_{n-1} X_{n-1} |^2}{|z_1|^2 + \cdots + |z_{n-1}|^2}$$ and the conclusion follows
by the Cauchy-Schwarz inequality.

Conversely, assume that $L(X, z)$ is  strictly positive
for every $(X_1,\dots ,X_{n-1}) \neq (0,\dots ,
0)$ and  each $z = (z_0,z_1,\dots ,z_{n-1})$
such that
$F(|z_0|^2) = |z_1|^2 +\cdots +
|z_{n-1}|^2$.
By inserting  $(X_1,\dots ,X_{n-1}) =
(z_1,\dots ,z_{n-1})$ in (\ref{Levi ristretto}) we get
$$L(z, z) = F(|z_0|^2) \left( 1 -
\frac{F' + F'' |z_0|^2}{F'^2 |z_0|^2} F(|z_0|^2) \right) >0$$
which implies $(xF' / F)' <0$.

\vskip 0.3cm

Finally, the proof of the equivalence
(ii)$\Leftrightarrow$(iv) is completely analogous to that given in
\cite{eng} (Proposition 3.4 and Proposition 3.6 ) for  the
$2$-dimensional case, to which the reader is referred. \fdim

\begin{remar}{\rm
Notice that the previous proposition  is a generalization   of
Proposition 3.6 in \cite{eng} proved there for the $2$-dimensional
case.}
\end{remar}

\vskip 0.3cm

Recall (see e.g.
\cite{kn}) that the Ricci curvature and the scalar curvature of a
\K\ metric $g$ on an $n$-dimensional complex manifold  $(M, g)$
are given
respectively by
\begin{equation}\label{riccicurvature}
{\Ric}_{\alpha \bar \beta} = - \frac{\partial^2}{\partial
z_{\alpha} \partial \bar z_{\beta}}(\log \det(h)), \ \ \alpha,
\beta = 0, \dots, n-1
\end{equation}
and
\begin{equation}\label{scalcurvature}
{\scal}_g= \sum_{\alpha, \beta = 0}^{n-1} g^{\beta \bar \alpha} {\Ric}_{\alpha \bar
\beta},
\end{equation}
where $g^{\beta \bar \alpha}$ are the entries  of the inverse
of  $(g_{\alpha\bar \beta})$, namely
$\sum_{\alpha =0}^{n-1}g^{\beta \bar \alpha}g_{\alpha\bar \gamma}=\delta_{\beta\gamma}$.

When $(M, g)=(D_F, g_F)$,
using  (\ref{matrixg})  it
is not
hard  to  check the validity of the following equalities.

\begin{equation}\label{g00}
g^{0 \bar 0} =\frac{A}{B} F,
\end{equation}

\begin{equation}\label{gb0}
g^{\beta \bar 0} =\frac{A}{B} F' z_0 \bar z_{\beta}, \ \ \ \
\beta = 1,\dots ,n-1,
\end{equation}

\begin{equation}\label{gba}
g^{\beta \bar \alpha} =\frac{A}{B} (F' + F'' |z_0|^2) z_{\alpha}
\bar z_{\beta}, \ \ \ \ \alpha \neq \beta, \ \alpha, \beta =
1,\dots ,n-1,
\end{equation}

\begin{equation}\label{gbb}
g^{\beta \bar \beta} =\frac{A}{B} [B + (F' + F'' |z_0|^2)
|z_{\beta}|^2], \ \ \ \ \beta = 1,\dots ,n-1,
\end{equation}

where
$$B = B(|z_0|^2)= F'^2 |z_0|^2 - F(F' + F'' |z_0|^2).$$
Now, set
$$L(x) = \frac{d}{dx} [x \frac{d}{dx} \log (xF'^2 -
F(F'+F''x))].$$
A straightforward computation using (\ref{determinante metrica})
and (\ref{riccicurvature}) gives:

\begin{equation}\label{ricci00}
{\Ric}_{0 \bar 0} = -L(|z_0|^2) - (n+1) g_{0 \bar 0},
\end{equation}

\begin{equation}\label{ricciAB}
{\Ric}_{\alpha \bar \beta} = - (n+1) g_{\alpha \bar \beta}, \ \ \
\alpha >0.
\end{equation}

Then, by (\ref{scalcurvature}),  the scalar curvature
of the metric $g_F$ equals
$${\scal}_{g_F} =- L(|z_0|^2) g^{0 \bar 0} -(n+1) \sum_{\alpha,\beta=0}^{n-1}
g^{\beta \bar\alpha} g_{\alpha \bar \beta} = - L(|z_0|^2) g^{0 \bar 0} -n(n+1) ,$$
which  by (\ref{g00}) reads as
\begin{equation}\label{scalgF}
{\scal}_{g_F}=  - \frac{A}{B} F L -n(n+1).
\end{equation}

\vskip 0.3cm

We conclude this section with 
Lemma \ref{mainlemma} below which will be used in the proof of our results.
This lemma is a
generalization of a result  proved by the first author
 for $2$-dimensional Hartogs domains
(see Theorem 4.8 in \cite{constscal}).
We first recall the definition of generalized scalar curvatures.
Given a  \K\ metric $g$ on   an  $n$-dimensional  complex manifold
$M$,  its {\em generalized scalar curvatures}
are the $n$ smooth functions
$\rho_0$,\dots ,$\rho_{n-1}$ on $M$
 satisfying  the following equation:
\begin{equation}\label{genscal}
\frac{\det(g_{\alpha \bar \beta} + t Ric_{\alpha \bar \beta})}{\det(g_{\alpha \bar \beta})} = 1 + \sum_{k=0}^{n-1} \rho_k t^{k+1} ,
\end{equation}
where $g_{\alpha\bar\beta}$ are the entries of the   metric in local coordinates.
Observe that for $k=0$ we recover the value of the scalar
curvature, namely
\begin{equation}\label{roscal}
\rho_0={\scal}_g.
\end{equation}
 The introduction and the study  of these curvatures  (in the compact case)
are due to K. Ogiue \cite{ogiue} to whom the reader is referred for
further results.  In particular,  in a joint paper with B.Y. Chen
\cite{CO}, he   studies  the constancy of one of the generalized
scalar curvatures. Their main result is that, under suitable
cohomological conditions,  the constancy of one of the
${\rho_k}'s, k=0, \dots ,n-1, $ implies that the metric $g$ is
Einstein.

\begin{lemma}\label{mainlemma}
Let $(D_F, g_F)$ be an $n$-dimensional  Hartogs domain. Assume
that one of its generalized scalar curvatures is constant. Then
$(D_F, g_F)$ is holomorphically isometric to an open subset of the  $n$-dimensional
hyperbolic space.
\end{lemma}
\dimostr
By  (\ref{ricci00}),
(\ref{ricciAB}) we get

$$\frac{det(g_{\alpha \bar \beta} + t Ric_{\alpha \bar
\beta})}{det(g_{\alpha \bar \beta})} = (1-(n+1)t)^{n} - t L
(1-(n+1)t)^{n-1} \frac{A F}{B} .$$

 So the generalized curvatures of $(D_F, g_F)$
 are given by

\begin{equation}
\rho_{k} = (n+1)^{k}(-1)^{k+1} {n-1 \choose k}\left[
\frac{n(n+1)}{k+1} + \frac{AFL}{B} \right], \ \ \ k = 0,\dots,n-1.
\end{equation}

Notice that,  for $k=0$,  we get $\rho_0= -\frac{AFL}{B}-n(n+1)=\scal_{g_F}$,
(compare with (\ref{scalgF})) in accordance with (\ref{roscal}).

Thus, $\rho_k$ is constant for some (equivalently, for any) $k
=0,\dots,n-1$ if and only if $\frac{AFL}{B}$ is constant. Since $A
= F(|z_0|^2) - |z_1|^2 - \cdots - |z_{n-1}|^2$ depends on
$z_1,\dots,z_{n-1}$ while $\frac{LF}{B}$ depends only on $z_0$,
this implies that $L =0$, i.e.

$$\frac{d}{dx} \left[x \frac{d}{dx} \log (xF'^2 - F(F'+F''x))\right]_{x = |z_0|^2} \equiv 0 .$$

Now, we continue as in the proof of Theorem 4.8 in
\cite{constscal} and conclude that $F (x)= c_1 - c_2 x, \ x=|z_0|^2$, with $c_1,
c_2 >0$, which implies that $D_F$ is holomorphically isometric
to an open subset of the complex hyperbolic space ${\C}H^{n}$ via the map

$$\phi: D_F \rightarrow {\C}H^{n}, \ (z_0, z_1,\dots,z_{n-1}) \mapsto \left( \frac{z_0}{\sqrt{c_1/c_2}}, \frac{z_1}{\sqrt{c_1}},\dots,\frac{z_{n-1}}{\sqrt{c_1}} \right) .$$
\fdim

\vskip 0.3cm

 \section{Proof of Theorem \ref{mainteor}}\label{final}

Extremal
metrics  were introduced and christened by Calabi \cite{caex} in
the compact case as the solution for the variational problem in  a
\K\ class defined by the square integral of the scalar curvature.
Therefore they are a generalization of constant scalar curvature
metrics. Calabi himself constructs nontrivial (namely with
nonconstant scalar curvature) metrics on some compact manifolds.
Only recently extremal \K\ metrics   were rediscovered  by several  mathematicians due to
their link with the  stability of complex vector bundles (see e.g.
 \cite{burns}, \cite{chtian},   \cite{fu}, \cite{levin} and  \cite{ma}).
Obviously extremal metrics  cannot be defined in the noncompact
case as  the solutions of a variational problem involving some
integral on the manifold. Nevertheless
they can be  alternatively defined (also in the noncompact case) as those  metrics
such that the (1,0)-part of the
Hamiltonian vector field associated to the scalar curvature is
holomorphic.
Therefore, in  local
coordinates an extremal metric must satisfy  
the following  system of  PDE's  (see \cite{caex}):
\begin{equation}\label{equazestrem}
\frac{\partial}{\partial \bar z_{\gamma}} \left( \sum_{\beta =
0}^{n-1} g^{\beta \bar \alpha} \frac{\partial \scal_g}{\partial
\bar z_{\beta}} \right) = 0,
\end{equation}
for every $\alpha, \gamma = 0,\dots,n-1$. 
Notice that in the noncompact case, the existence and
uniqueness of  such    metrics are far from being  understood. For
example, only recently in \cite{cheng} (see also \cite{cheng2}),
 it has been shown the existence of a nontrivial
extremal   and complete  \K\ metric
 in a complex one-dimensional manifold.

\vskip 0.3cm

\noindent
{\bf Proof of   Theorem \ref{mainteor}}
 \noindent
 In order to use  equations (\ref{equazestrem}) for $(D_F, g_F)$
we write   the entries  $g^{\beta \bar \alpha}$
by separating the terms depending  only on  $z_0$ from the other terms.   More precisely, (\ref{g00}), (\ref{gb0}), (\ref{gba}) and (\ref{gbb})  can be written as follows.

$$g^{0 \bar 0} = P_{00} + Q_{00}(|z_1|^2+ \cdots +|z_{n-1}|^2), $$

$$g^{0 \bar \alpha} = \bar z_0 z_{\alpha} [P_{0a} + Q_{0a}(|z_1|^2+\cdots +|z_{n-1}|^2)], \ \ \alpha=1,\dots ,n-1, $$

$$g^{\alpha \bar \alpha} = F + P_{aa}|z_{\alpha}|^2 - (1 + Q_{aa}|z_{\alpha}|^2)\sum_{k \neq \alpha} |z_k|^2 - R_{aa}|z_{\alpha}|^4, \ \ \alpha=1,\dots ,n-1, $$

$$g^{\beta \bar \alpha} = \bar z_{\beta} z_{\alpha}[P_{ab} + Q_{ab}(|z_1|^2+\cdots +|z_{n-1}|^2)],\ \alpha \neq \beta, \ \alpha, \beta=1,\dots ,n-1, $$

where

\label{definizioniPQR}

$$P_{00} = \frac{F^2}{B}, \ \ \ Q_{00} = - \frac{F}{B}, $$

$$P_{0a} = \frac{F' F}{B}, \ \ \ Q_{0a} = - \frac{F'}{B},  $$

$$P_{aa} = \frac{F(F' + F'' |z_0|^2)}{B} -1 , \ \ \ Q_{aa} = R_{aa} =  \frac{F' + F'' |z_0|^2}{B} , $$

$$P_{ab} = \frac{F(F' + F'' |z_0|^2)}{B} , \ \ \ Q_{ab} =  - \frac{F' + F'' |z_0|^2}{B} $$

are all functions depending only on $|z_0|^2$.

We also have (cfr. (\ref{scalgF}))

\begin{equation}\label{scalG}
{\scal}_{g_F} = -n(n+1) + G (F - |z_1|^2 -\cdots -|z_{n-1}|^2)
\end{equation}

where
$$G=G(|z_0|^2) = - \frac{L(|z_0|^2) F(|z_0|^2)}{B(|z_0|^2)}.$$

\vskip 0.3cm

Assume that $g_F$ is an extremal metric, namely equation
(\ref{equazestrem}) is satisfied. We are going to show that
$\scal_{g_F}$ is constant and hence by Lemma \ref{mainlemma}
$(D_F, g_F)$ is holomorphically isometric to an open subset of  $({\complex}H^{n},
g_{hyp})$. In order to do that, fix  $i \geq 1$ and let us  consider
equation (\ref{equazestrem}) when $g=g_F$ for $\alpha = 0$,
$\gamma = i$.

We have

$$\frac{\partial \scal_{g_F}}{\partial \bar z_0} = G' z_0 (F - |z_1|^2 -\cdots -|z_{n-1}|^2) + z_0 G F'$$

$$\frac{\partial \scal_{g_F}}{\partial \bar z_i} = -G z_i. $$

So, equation (\ref{equazestrem}) gives

$$\frac{\partial}{\partial \bar z_i}\left\{ \left[P_{00} + Q_{00}\sum_{k=1}^{n-1} |z_k|^2 \right]\left[ G' z_0 (F - \sum_{k=1}^{n-1} |z_k|^2) + z_0 G F' \right] \right. -$$
$$\left.- z_0 G \left[P_{0a} + Q_{0a}\sum_{k=1}^{n-1} |z_k|^2\right]
\sum_{k=1}^{n-1} |z_k|^2 \right\} = 0,$$

namely

$$Q_{00}z_i \left[ G' z_0 (F - \sum_{k=1}^{n-1} |z_k|^2) + z_0 G F' \right] - G' z_0 z_i \left[P_{00} + Q_{00}\sum_{k=1}^{n-1} |z_k|^2 \right] -$$

$$ -z_0 G Q_{0a} z_i \sum_{k=1}^{n-1} |z_k|^2 - z_0 z_i G \left[P_{0a} + Q_{0a}\sum_{k=1}^{n-1} |z_k|^2\right] = 0 .$$

Deriving again  with respect to $\bar z_i$,  we get
$$ -2
Q_{00} G' z_0 z_i^2 - 2 G Q_{0a} z_0 z_i^2 = 0.$$ Let us assume
$z_0z_i \neq 0$. This  implies $Q_{00}
G' + G Q_{0a} = 0$, i.e.  $GF' + FG' = 0$ or, equivalently,
 $G = \frac{c}{F}$ for some constant $c \in {\R}$. The
proof of Theorem \ref{mainteor} will be completed by showing that
$c =0$.  In fact, in this case $G =0$ on the open and dense subset
of $D_F$ consisting of those points such that $z_0z_i \neq 0$
and therefore, by (\ref{scalG}),
$\scal_{g_F}$ is   constant on  $D_F$.
 In order to prove that $c=0$, let us now consider
equation (\ref{equazestrem}) for $\alpha = i$, $\gamma = i$.

$$\frac{\partial}{\partial \bar z_i} \left\{ \bar z_0 z_i \left[ G' z_0
( F - \sum_{k=1}^{n-1} |z_k|^2 ) + G F' z_0 \right] \left[P_{0a} +
Q_{0a}\sum_{k=1}^{n-1} |z_k|^2  \right]  - \right.$$

$$ - G z_i \left[ F + P_{aa}|z_i|^2 -
(1 + Q_{aa}|z_i|^2)\sum_{k \neq 0,i} |z_k|^2 - R_{aa}|z_i|^4
\right] -$$
$$\left.- G z_i \sum_{k \neq 0, i} |z_k|^2 \left[ P_{ab} +
Q_{ab}\sum_{k=1}^{n-1} |z_k|^2 \right] \right\}  =0 .$$

This implies

$$ -G' |z_0|^2 z_i^2 \left[P_{0a} + Q_{0a}\sum_{k=1}^{n-1} |z_k|^2 \right] + \bar z_0 z_i^2 Q_{0a} \left[ G' z_0 ( F - \sum_{k=1}^{n-1} |z_k|^2 ) + G F' z_0 \right] -  $$

$$- P_{aa} G z_i^2 + G z_i^2 Q_{aa}  \sum_{k \neq 0, i} |z_k|^2 + 2 G z_i^3 \bar z_i R_{aa} - G z_i^2 Q_{ab}  \sum_{k \neq 0, i} |z_k|^2. $$

If we divide by $z_i^2$ (we are assuming $z_i \neq 0$) and derive
again the above expression with respect to $\bar z_i$ we get
$$ -G' |z_0|^2 Q_{0a} + G R_{aa} = 0 .$$
By the definitions made at page \pageref{definizioniPQR}
this is equivalent to

$$ \frac{G' F' |z_0|^2 + G (F' + F'' |z_0|^2)}{B} = 0, $$
i.e. $(G F' x)' = 0, x=|z_0|^2$. Substituting $G = \frac{c}{F}$
in this equality we get $c (\frac{F' x}{F})' = 0$.
Since $(\frac{F' x}{F})' <0$
(by (ii) in Proposition \ref{Kmet}) $c$ is forced to be zero,
and this concludes the proof.
\fdim

\section{Proof of Theorem \ref{mainteor2}}\label{final2}
A {\em \K--Ricci soliton} on a complex manifold $M$ is a pair $(g,
X)$ consisting of a \K\ metric $g$ and a real holomorphic vector
field $X$ on $M$ such that
\begin{equation}\label{eqsol}
{\Ric}_g = \lambda \ g + L_Xg,
\end{equation}
for some $\lambda \in \R$, where $L_X g$ is the Lie derivative of
$g$ along $X$, i.e.
\begin{equation}\label{lieder}
(L_Xg)(Y,Z) = X(g(Y,Z)) - g([X,Y], Z) - g(Y, [X,Z]),
\end{equation}
for $Y, Z$ vector fields  on $M$.  A {\em real holomorphic vector
field} $X$  is the real part of a holomorphic vector field,
namely,
 in  local complex coordinates $(z_0,\dots,z_{n-1})$ on an open subset  $U\subset M$,
\begin{equation}\label{holvectfield}
X = \sum_{k=0}^{n-1} \left( f_k \frac{\partial}{\partial z_k} + \bar
f_k \frac{\partial}{\partial \bar z_k} \right),
\end{equation}
for some holomorphic functions $f_k$, $k=0,\dots,n-1$ on $U$.

We refer the reader to \cite{tiansol1}, \cite{tiansol2},
\cite{tiansol3} for the existence and uniqueness of \K--Ricci
solitons on compact manifolds and to \cite{FI} for the noncompact
case. \K--Ricci solitons generalize \K--Einstein metrics. Indeed any
\K--Einstein metric $g$ on a complex manifold $M$ gives rise to a
trivial \K--Ricci soliton by choosing $X = 0$ or $X$ Killing with
respect to $g$. Obviously if the automorphism group of $M$ is
discrete then a \K--Ricci soliton $(g, X)$ is nothing but a
\K--Einstein metric $g$.

Our Theorem \ref{mainteor2} asserts that a \K--Ricci soliton
$(g_F, X)$ on a pseudoconvex Hartogs domain $D_F$ is necessarily
trivial. Notice that the automorphism group of $D_F$ is not
discrete (see also \cite{iskra}).

\vskip 0.3cm

\noindent
{\bf Proof of   Theorem \ref{mainteor2}}
 \noindent
 Let $(g_F, X)$ be a \K--Ricci soliton.
 By applying both sides  of  (\ref{eqsol}) to the pair $(\frac{\partial}{\partial z_0},\frac{\partial}{\partial\bar z_0})$
 and taking into account  (\ref{ricci00}) one gets:
\begin{equation}
-L(|z_0|^2) = \gamma g_{0 \bar 0} + \sum_{k=0}^{n-1} \left( f_k
\frac{\partial g_{0 \bar 0}}{\partial z_k} + \bar f_k
\frac{\partial g_{0 \bar 0}}{\partial \bar z_k} \right) +
\sum_{k=0}^{n-1} \left( \frac{\partial f_k}{\partial z_0} g_{k
\bar 0} + \frac{\partial \bar f_k}{\partial \bar z_0} g_{0 \bar k}
\right)
\end{equation}

\noindent where

\begin{equation}\label{gamma}
\gamma = \lambda + (n+1).
\end{equation}

By (\ref{matrixg}), we have

\begin{equation}\label{equation2}
\tilde C = \sum_{k=0}^{n-1} C_k (f_k \bar z_k + \bar f_k z_k) +
C(\phi_0 + \bar \phi_0) - F' \sum_{k=1}^{n-1} \left( z_0 \bar z_k
\frac{\partial f_k}{\partial z_0} +\bar z_0 z_k \frac{\partial
\bar f_k}{\partial \bar z_0} \right)
\end{equation}

\noindent \noindent where we have set $\tilde C = -A^2 L - \gamma
C$, $C_k = A^2 \frac{\partial g_{0 \bar 0}}{\partial x_k}$ ($x_k =
|z_k|^2$) and $\phi_0 = \frac{\partial f_0}{\partial z_0}$ ($A$
and $C$ are given by (\ref{eqA}) and (\ref{eqc}) respectively).

Now, by applying the operator
$\frac{\partial^4}{\partial^2 z_i
\partial^2 \bar z_i}$ ($i = 1, \dots, n-1$) to both sides of this
equation we get
\begin{small}

\begin{eqnarray}
-4 L & = & \sum_{k=0}^{n-1} \frac{\partial^4 C_k}{\partial z_i^2
\partial \bar{z_i}^2}(f_k \bar{z_k} + \bar{f_k} z_k) +
2\sum_{k=0}^{n-1} \frac{\partial^3 C_k}{\partial z_i^2
\partial \bar{z_i}}\left( f_k \delta_{ik} + \frac{\partial \bar{f_k}}{\partial \bar{z_i}} z_k
\right)+ \nonumber \\ & + &  2\sum_{k=0}^{n-1} \frac{\partial^3
C_k}{\partial z_i
\partial \bar{z_i}^2}\left( \bar{f_k} \delta_{ik} + \frac{\partial f_k}{\partial z_i} \bar{z_k} \right) + 4 \sum_{k=0}^{n-1} \frac{\partial^2 C_k}{\partial z_i
\partial \bar{z_i}}\left( \frac{\partial f_k}{\partial z_i} \delta_{ik} +  \frac{\partial \bar{f_k}}{\partial \bar{z_i}} \delta_{ik}
\right) + \nonumber \\ & + & \sum_{k=0}^{n-1} \frac{\partial^2
C_k}{\partial z_i^2} \frac{\partial^2 \bar{f_k}}{\partial
\bar{z_i}^2} z_k + \sum_{k=0}^{n-1} \frac{\partial^2 C_k}{\partial
\bar{z_i}^2} \frac{\partial^2 f_k}{\partial z_i^2} \bar{z_k} +
2\sum_{k=0}^{n-1} \frac{\partial C_k}{\partial z_i}
\frac{\partial^2 \bar{f_k}}{\partial \bar{z_i}^2} \delta_{ik} + \nonumber \\
& + & 2\sum_{k=0}^{n-1} \frac{\partial C_k}{\partial \bar{z_i}}
\frac{\partial^2 f_k}{\partial z_i^2} \delta_{ik} +
\frac{\partial^4 C}{\partial z_i^2 \partial \bar{z_i}^2}(\phi_0 +
\bar{\phi_0}) + 2\frac{\partial^3 C}{\partial z_i^2 \partial
\bar{z_i}} \frac{\partial \bar{\phi_0}}{\partial \bar{z_i}} + \nonumber \\
& + & 2\frac{\partial^3 C}{\partial z_i \partial \bar{z_i}^2}
\frac{\partial \phi_0}{\partial z_i} + \frac{\partial^2
C}{\partial \bar{z_i}^2} \frac{\partial^2 \phi_0}{\partial z_i^2}+
\frac{\partial^2 C}{\partial z_i^2} \frac{\partial^2
\bar{\phi_0}}{\partial \bar{z_i}^2} \nonumber.
\end{eqnarray}

\end{small}

\noindent Since $C$ and $C_k$ are rotation invariant, by
evaluating the previous expression at $z_1= \cdots = z_{n-1} = 0$
and taking into account that

$$\frac{\partial^4
C_0}{\partial z_i^2 \partial
\bar{z_i}^2}|_{\{ z_1= \cdots = z_{n-1} = 0 \}} = - 8
x \frac{F'^3}{F^3},$$

$$\frac{\partial^2 C_i}{\partial z_i \partial
\bar{z_i}}|_{\{ z_1= \cdots = z_{n-1} = 0 \}} = 2 x
\frac{F'^2}{F^2},$$

$$\frac{\partial^4 C}{\partial z_i^2 \partial
\bar{z_i}^2}|_{\{ z_1= \cdots = z_{n-1} = 0 \}} = 0,$$

we have

\begin{equation}\label{00derivato}
L = 2 x \frac{F'^3}{F^3} (f_0
\bar{z_0} +
\bar{f_0} z_0) -  2 x \frac{F'^2}{F^2} (\phi_i +
\bar{\phi}_i),
\end{equation}

where $\phi_ i = \frac{\partial f_i}{\partial z_i}$.

\medskip

\noindent Now, let $i = 1, \dots, n-1$.  By applying both sides of
(\ref{eqsol}) to the pair $(\frac{\partial}{\partial
z_{i}},\frac{\partial}{\partial\bar z_{i}})$ one gets

\begin{equation}\label{alfaalfa}
-\gamma \ g_{i \bar i} = \sum_{k=0}^{n-1} \left( f_k \frac{\partial
g_{i \bar i}}{\partial z_k} + \bar f_k \frac{\partial g_{i \bar
i}}{\partial \bar z_k} \right) + \sum_{k=0}^{n-1} \left(
\frac{\partial f_k}{\partial z_{i}} g_{k \bar {i}} +
\frac{\partial \bar f_k}{\partial \bar z_{i}} g_{{i} \bar k}
\right)
\end{equation}

\noindent where $\gamma$ is given by (\ref{gamma}). By
(\ref{matrixg}) and (\ref{ricciAB}) this means

\begin{small}

\begin{eqnarray}\label{alfaalfa2}
-\gamma \ \frac{A + |z_{i}|^2}{A^2} & = & \sum_{k=1}^{n-1} A^{-3}
[2 |z_{i}|^2 + A(1 + \delta_{i k})] ( f_k \bar z_k + \bar f_k z_k)
- \nonumber \\ & - & A^{-3} F' (2 |z_{i}|^2 + A)(f_0 \bar z_0 +
\bar f_0 z_0) - \frac{F'}{A^2} \left( \frac{\partial f_0}{\partial
z_{i}} \bar z_0 z_{i} + \frac{\partial \bar f_0}{\partial \bar
z_{i}} z_0 \bar z_{i} \right) + \nonumber \\ & + & \frac{1}{A^2}
\sum_{k=1}^{n-1} \left[ \frac{\partial f_k}{\partial z_{i}} (\bar
z_k z_{i} + \delta_{k i} A) + \frac{\partial \bar f_k}{\partial
\bar z_{i}} (z_k \bar z_{i} + \delta_{k i} A) \right].
\end{eqnarray}

\end{small}

If we evaluate both sides of this equation at $z_1 =
\cdots = z_{n-1} = 0$ we get

\begin{equation}\label{equgamma}
-\gamma F = - F' (f_0 \bar z_0 + \bar f_0 z_0) + F (\phi_{i} +
\bar \phi_{i}).
\end{equation}

Moreover,   by multiplying equation (\ref{alfaalfa2}) by
$A^2$, by  applying the operator $\frac{\partial^2}{\partial z_{i}
\partial \bar z_{i}}$ to both sides and evaluating  at $z_1 =
\cdots = z_{n-1} = 0$ one gets
\begin{equation}\label{last}
0 = - \frac{F'}{F} (f_0 \bar z_0 + \bar f_0 z_0) + (\phi_{i} +
\bar \phi_{i}).
\end{equation}

Finally, by comparing  (\ref{00derivato}) with  (\ref{last}),
one gets   $L=0$ and hence, by the proof of Lemma \ref{mainlemma},
 $(D_F, g_F)$ is holomorphically
isometric to an open subset of $({\C} H^n, g_{hyp})$ and we are done.
(Notice that  equations    (\ref{equgamma}) and (\ref{last})  yield $\gamma =0$ and  by (\ref{eqsol}) with $g_F=g_{hyp}$  one  gets that $X$ is a  Killing vector field with respect to $g_{hyp}$, as expected).
\fdim

\small{}

\end{document}